\documentclass[12pt]{elsarticle}
\usepackage{amssymb}
\usepackage{amsfonts}
\usepackage{color,graphicx}
\usepackage{latexsym}
\usepackage{amssymb}
\newtheorem{thm}{Theorem}[section]

\newtheorem{?}{?}
\newtheorem{coro}[thm]{Corollary}

\newcommand{\Z}{\mathbb{Z}}

\newcommand{\beeq}{\begin{eqnarray*}}
\newcommand{\eneq}{\end{eqnarray*}}
\newcommand{\proof}{\noindent {\it Proof.\hspace{4mm}}}

\newcommand{\qfd}{\hfill $\fbox{}$\vspace{4mm}}\def\newpic#1{%
\def\emline##1##2##3##4##5##6{%
\put(##1,##2){\special{em:point #1##3}}%
\put(##4,##5){\special{em:point #1##6}}%
\special{em:line #1##3,#1##6}}}
\newpic{}
\def\emline#1#2#3#4#5#6{%
\put(#1,#2){\special{em:moveto}}%
\put(#4,#5){\special{em:lineto}}}
\def\newpic#1{}
\usepackage{lineno, hyperref}
\begin{document}
\begin{frontmatter}
\title{On Coloring the Arcs of Biregular Graphs}
\author{Italo J. Dejter}

\address{University of Puerto Rico, Rio Piedras, Puerto Rico 00936--8377}

\ead{italo.dejter@gmail.com}

\begin{abstract}
%\begin{linenumbers}
\noindent
Recalling each edge of a graph $H$ has 2 oppositely oriented arcs, each vertex $v$ of $H$ is identified with the set of arcs, denoted $(v,e)$, departing from $v$ along the edges $e$ of $H$ incident to $v$. Let $H$ be a $(\lambda,\mu)$-biregular graph with bipartition $(Y,X)$, where $|Y|=k\mu$ and $|X|=k\lambda$, ($0<k,\lambda,\mu\in\Z$). We consider the problem, for each edge $e=yx$ in $H$, of assigning, a color (given by an element) of $Y$, resp. $X$, to the  arc $(y,e)$, resp. $(x,e)$, so that each color is assigned exactly once in the set of arcs departing from each vertex of $H$. Furthermore, we set such assignment to fulfill a specific bicolor weight function over a monotonic subset of $Y\times X$. This problem applies to the Design of Experiments for Industrial Chemistry, Molecular Biology, Cellular Neuroscience, etc. An algorithmic construction based on biregulzr graphs with bipartitions given by cyclic-group pairs is presented, as well as 3 essentially different solutions to the  Great Circle Challenge Puzzle based on a different biregular graph whose bipartition is formed by the vertices and 5-cycles of the Petersen graph.
%\end{linenumbers}
\end{abstract}

\end{frontmatter}

\section{Introduction}\label{s1}

Assume we need to perform a sequence of Industrial Chemistry experiments involving a model represented by a bipartite graph $H$ with parts $Y$ and $X$ and in which each $y\in Y$ must react with each of its neighbors $x\in X$ and vice versa. Say these experiments involve a number of {\it acids} equal to the degree of each $y$ and a number of {\it bases} equal to the degree of each $x$. Moreover, assume an ordered presentation of acid-base pairs together with a weight function (representing, say, the number of different acid-base required proportions). Our task would be to distribute the experiments so that the pairs formed by a member of $Y$ and a member of $X$ contribute to fulfill the weight function requirements. 

This combinatorial problem can be interpreted as an assemblage one, like in the concrete application to the Great Circle Challenge Puzzle presented in Section~\ref{s5} or in the Design of Experiments in Industrial Chemistry (as suggested above), Molecular Biology, Cellular Neuroscience, etc.

In a Cellular Neuroscience setting, the vertices of $Y$ and $X$ may be thought of as representing {\it source} and {\it target} neurons respectively, while the arcs from $Y$ to $X$ may be assigned {\it transmitter proteins} as their colors, and the arcs from $X$ to $Y$ may be assigned {\it receptor proteins} as their colors. The edges between $Y$ and $X$ themselves may be thought of as {\it synapses}. 

Let us formalize these ideas.
If $0<k,\lambda,\mu\in\mathbb{Z}$ with $\lambda\ne\mu$, then a bipartite graph $H=(Y\cup X,E_H)$ with bipartition $(Y,X)$ is said to be $(\lambda,\mu)$-{\it biregular} (\cite{SU}, page 137) or $(\lambda,\mu)$-{\it semiregular}  (\cite{GR}, page 12) if every $y\in Y$, resp. $x\in X$, has degree $\lambda$, resp. $\mu$, where $|X|=k\lambda$ and $|Y|=k\mu$.

We will be working with a coloring {\it palette} $\Theta\subseteq\Z_\lambda \times\Z_\mu  $, where $\Z_\lambda =\{0,\ldots,\lambda-1\}$, resp. $\Z_\mu  =\{0,\ldots,\mu-1\}$, is taken as a {\it color set} for $Y$, resp. $X$, (representing, say,  acids or transmitter proteins, resp. bases or receptor proteins), with $0<|\Theta|<\lambda|Y|=\mu|X|$.

Given $e=yx\in E_H$ with $(y,x)\in Y\times X$, the pairs $(y,e)$ and $(x,e)$ are said to be the {\it arcs} of $e$, with {\it tails} $y$ and $x$ and {\it heads} $x$ and $y$, respectively.

 A $(\lambda,\mu)$-{\it coloring} of $H$ by $\Theta$ is an assignment of a {\it bicolor} $(\alpha_y^e,\beta_x^e)\in\Theta$ to each $e=yx\in E_H$, meaning both that $\alpha_y^e$ is assigned to the arc $(y,e)$, and that $\beta_x^e$ is assigned to the arc $(x,e)$. We say that such a coloring is {\it full} if the following two conditions are satisfied:

\noindent{\bf 1.} $\forall(y,\alpha)\in Y\times\Z_\lambda , \exists e=yx\in E_H$ with assigned bicolor $(\alpha,\beta_x^e)\in\Theta$,  for a

sole $x\in X$; (i.e. $\alpha$ is assigned to $(y,e)$, and also $\exists\beta_x^e$ assigned to $(x,e)$);

\noindent{\bf 2.} $\forall(x,\beta)\in X\times\Z_\mu  , \exists e=yx\in E_H$ with assigned bicolor $(\alpha_y^e,\beta)\in\Theta$, for a

sole $y\in Y$; (i.e. $\beta$ is assigned to $(x,e)$, and also $\exists\alpha_y^e$ assigned to $(y,e)$).

\noindent
In all full $(\lambda,\mu)$-colorings considered below, $\Theta$ also satisfies the {\it monotonicity} property, namely:
if $(\alpha,\beta)$ and $(\alpha',\beta')$ are in $\Theta$, then $\alpha\le \alpha'$ if and only if $\beta\le\beta'$.

Let us consider the problem of determining full $(\lambda,\mu)$-colorings in specific $(\lambda,\mu)$-biregular graphs.
Each such coloring determines a bijective correspondence from the set of arcs departing from each vertex of $Y$, resp. $X$, onto $\Z_\lambda$, resp. $\Z_\mu$.
Given one such $(\lambda,\mu)$-coloring, we define the {\it weight function} 

\begin{equation}\omega:\Theta\rightarrow\{1,\ldots,|E_H|\}\end{equation}

\noindent to express the number $\omega(\alpha,\beta)$ of edges with bicolor $(\alpha,\beta)\in\Theta$ and restrict the problem by imposing the natural constraint 
 
 $$\sum_{(\alpha,\beta)\in\Theta}\omega(\alpha,\beta)=|E_H|.$$

In Section~\ref{s2}, this problem leads, via Theorem~\ref{thm1} and Corollaries~\ref{co1}-\ref{co2}, to an algorithmic characterization of full $(\lambda,\mu)$-colorings of $(\lambda,\mu)$-biregular graphs based on cyclic-group pairs, for $\lambda<\mu$ and $\lambda\nmid\mu$. 
In Corollary~\ref{co2}, a $d$-plicating extension (duplicating extension if $d=2$, exemplified from the penultime to the last $(\lambda,\mu)$-list in (6)) shows that in the realm of cyclic-group pairs the problem is solved.
Yet, the decision problem of determining whether there exists, or not, a full $(\lambda,\mu)$-coloring of any particular instance of a $(\lambda,\mu)$-biregular graph is wide open and possibly hard, as is the associated counting problem, considering the connections between the mentioned applicative problems and purely mathematical aspects.
Away from the approach via cyclic-group pairs in Section~\ref{s2}, and in order to augment involvement with decision and counting,
odd-girth cycles in odd graphs \cite{B,BCN,GR} (see Sections~\ref{s5/2}-\ref{s3} below) are shown to lead
to a $(5,3)$-biregular graph $H'$ based on the vertices and 5-cycles of the Petersen graph $O_3$ in which we counted three full $(5,3)$-colorings, of which only one is reversible (Theorem~\ref{abc}).

The last $(5,3)$-biregular graph $H$ cited in (6) has one such $(5,3)$-coloring.
But $H'$ and $H$ are nonisomorphic. 
In fact, the automorphism groups of all graphs in Section~\ref{s2} are abelian, while that of $H'$ is nonabelian, making it interesting.
So
in Sections~\ref{s4}-\ref{s5}, it is shown that the bipartite double cover $G'$ of $H'$ (obtained by duplicating extension from $H'$) yields a model of the Great Circle Challenge Puzzle (or GCCP) assemblage problem (available online) conceived by Buckminster Fuller's disciple Tim Tendick.
As commented above, three solutions to this GCCP assemblage problem arise (in Corollary~\ref{coro}, from Theorems~\ref{abc}-\ref{assem}), only one being reversible.

\section{Approach via Cyclic-Group Pairs}\label{s2}

Let $m,n\in\Z$ with $n<m$ and $n\nmid m$. The notation in the following statement is adopted: 
Each $i\in\Z$ determines unique elements 

\begin{equation}i_m=(i\mbox{ mod }m)\in\Z_m\;\mbox{ and }\;i_n=(i\mbox{ mod }n)\in\Z_n.\end{equation} 

Let $H$ be a biregular graph with bipartition $(Y,X)=(\Z_m,\Z_n)$ and an edge 

\begin{equation}i_mi_n=(i\mbox{ mod }m)(i\mbox{ mod }n)\in E_H\end{equation}

\noindent for each $i\in\{0,\ldots,\frac{mn}{\gcd(m,n)}-1\}$.
It is easy to see that the regular degrees of $Y$ and $X$ in $H$ are $\lambda=\frac{n}{\gcd(m,n)}$ and $\mu=\frac{m}{\gcd(m,n)}$, respectively. 
To see that $H$ is connected, notice that if we assumed that $n|m$ (contrary to our initial assumption $n\nmid m$), then $H$ would be disconnected, specifically with exactly $\frac{m}{n}$ components, where $\frac{m}{n}>1$. (These components are induced by the disjoint vertex sets $\{(i+jn)_m;0\le j<\frac{m}{n}\}$, for $0\le i<n$).

As in (2), each $i\in\Z$ determines unique elements $i_\lambda=(i$ mod $\lambda)\in\Z_\lambda$ and $i_\mu=(i$ mod $\mu)\in\Z_\mu$. 
For the instance $(m,n)=(9,6)$, where $\lambda=2$ and $\mu=3$, the assignment (5), below, will be a (2,3)-coloring of $H$ by the set

\begin{equation}\Theta=\{(0_2,0_3),(0_2,1_3),(1_2,1_3),(1_2,2_3)\}\subset\Z_2\times\Z_3,\end{equation} 

\noindent having (see display (1)) the weight function $\omega:\Theta\rightarrow\{1,\ldots,18\}$ given by  $\omega(0_2,0_3)=6,\omega(0_2,1_3)=3,\omega(1_2,1_3)=3$, $\omega(1_2,2_3)=6$. But by eliminating the subindices 2 and 3  and just writing $\Theta=\{(0,0),(0,1),(1,1),(1,2)\}$ instead of the expression (4), the assignment in question can be written as:

\begin{equation}\begin{array}{llllll}
^{0_90_6\rightarrow(0,0),}_{6_90_6\rightarrow(0,1),}&^{1_91_6\rightarrow(0,0),}_{7_91_6
\rightarrow(0,1),}&^{2_92_6\rightarrow(0,0),}_{8_92_6\rightarrow(0,1),}&^{3_93_6
\rightarrow(0,0),}_{0_93_6\rightarrow(1,1),}&^{4_94_6\rightarrow(0,0),}_{1_94_6
\rightarrow(1,1),}&^{5_95_6\rightarrow(0,0),}_{2_95_6\rightarrow(1,1),}\\
^{3_90_6\rightarrow(1,2),}&^{4_91_6\rightarrow(1,2),}&^{5_92_6\rightarrow(1,2),}
&^{6_93_6\rightarrow(1,2),}&^{7_94_6\rightarrow(1,2),}&^{8_95_6\rightarrow(1,2).}
\end{array}\end{equation}

This assignment is {\it reversible}, meaning that the binary actions on $\Z_\lambda $ and $\Z_\mu  $ given by the permutations sending $\Z_\lambda $ and $\Z_\mu  $ respectively onto $(\lambda-1,\lambda-2,\ldots,1)$ and $(\mu-1,\mu-2,\ldots,1)$ still yield a $(\lambda,\mu)$-coloring of $H$.

\vspace*{2mm}

\begin{thm}\label{thm1} Given positive integers $\lambda,\mu$ and $k$ with $\lambda<\mu$ and $\lambda\nmid\mu$, let
$(m,n)=(k\mu,k\lambda)$. Then,
the biregular graph $H$ with bipartition $(Y,X)=(\mathbb{Z}_m,\mathbb{Z}_n)$ and edges expressible as in (3) admits a reversible full $(\lambda,\mu)$-coloring by $\Theta=\{(0,0),(0,1),\ldots,(\mu-1,\lambda-1)\}$. 
\end{thm}

\proof An ordered list 
$L=\{0_m0_n,1_m1_n,\ldots,i_mi_ n,\ldots,(m-1)_m(n-1)_n\}$ can be readily built that contains each edge of $H$ exactly once. This allows to construct a full $(\lambda,\mu)$-coloring for $H$ via de following algorithm:

\noindent{\bf(1)} assign bicolor $(0_\lambda,0_\mu)$ to the first edge $i_mi_n=0_m0_n$ of $L$;

\noindent{\bf(2)} if edge $i_mi_n$ in $L$ was assigned bicolor $(\alpha_\lambda,\beta_\mu)$, then assign to the following

edge $(i+1)_m(i+1)_n$ of $L$ the following bicolor:

{\bf(a)} $(\alpha_\lambda,\beta_\mu)$, if no contradiction to items 1.-2. in Section~\ref{s1}, occurs; else:

{\bf(b)} either($\alpha_\lambda+1_\lambda,\beta_\mu)$ or $(\alpha_\lambda,\beta_\mu+1_\mu)$, where necessarily only one of $\alpha_\lambda$

\hspace*{8mm} and $\beta_\mu$ does not contradict items 1.-2. in Section~\ref{s1};

\noindent{\bf(3)} let $i+1$ take the place of $i$ and go to step (3), unless $i=\frac{mn}{\gcd(m,n)}$, in

which case the algorithm stops, yielding a set  $\Theta$ for which a monotonic 

full $(\lambda,\mu)$-coloring of $H$ can be reconstructed.

\noindent It is easy to check that this produces a reversible $(\lambda,\mu)$-coloring by $\Theta=\{(0_\lambda,0_\mu),(0_\lambda,1_\mu),\ldots,(\mu_\lambda-1_\lambda,\lambda_\mu-1_\mu)\}$, as exemplified subsequently.
\qfd

\noindent Examples of the running of the algorithm in the proof of Theorem~\ref{thm1} are shown
in the following three tables for $(m,n)=(18,12),(15,6),(15,12)$,
where: {\bf(i)} the bicolors $(\alpha_\lambda,\beta_\mu)$, shortened as $\alpha\beta$, are shown on the leftmost column in each table; {\bf(ii)} the endvertex $i_n$ of each edge $i_mi_n$ in $L$ heads respectively the corresponding second, third, $\ldots$, last column; {\bf(iii)} the endvertex $i_m$ of each edge $i_mi_n$  in $L$ is located below the column header $i_n$ on the row headed on the left by the corresponding assigned bicolor $\alpha\beta$; {\bf(iv)} hexadecimal notation is used (e.g. $a=10,b=11$):

$$\begin{array}{||c|cccccc||c|cccccc||c|cccccccccccc||}
^{\alpha\beta}_{==}\!\!&\!^{0_6}_{=}\!\!&\!^{1_6}_{=}\!\!&\!^{2_6}_{=}\!\!&\!^{3_6}_{=}\!\!&\!^{4_6}_{=}\!\!&\!^{5_6}_{=}\!\!&\!^{\alpha\beta}_{==}\!\!&\!^{0_6}_{=}\!\!&\!^{1_6}_{=}\!\!&\!^{2_6}_{=}\!\!&\!^{3_6}_{=}\!\!&\!^{4_6}_{=}\!\!&\!^{5_6}_{=}
\!&\!\!^{\alpha\beta}_{==}\!\!&\!^{0_c}_{=}\!\!&\!^{1_c}_{=}\!\!&\!^{2_c}_{=}\!\!&\!^{3_c}_{=}\!\!&\!^{4_c}_{=}\!\!&\!^{5_c}_{=}\!\!&\!^{6_c}_{=}\!\!&\!^{7_c}_{=}\!\!&\!^{8_c}_{=}\!\!&\!^{9_c}_{=}\!\!&\!^{a_c}_{=}\!\!&\!^{b_c}_{=}\\
^{00}_{01}\!\!&\!^{0_9}_{6_9}\!\!&\!^{1_9}_{7_9}\!\!&\!^{2_9}_{8_9}\!\!&\!^{3_9}_{-}\!\!&\!^{4_9}_{-}\!\!&\!^{5_9}_{-}\!\!
&\!\!^{00}_{01}\!\!&\!^{0_f}_{6_f}\!\!&\!^{1_f}_{7_f}\!\!&\!^{2_f}_{8_f}\!\!&\!^{3_f}_{9_f}\!\!&\!^{4_f}_{a_f}\!\!&\!^{5_f}_{b_f}
\!&\!\!^{00}_{01}\!\!&\!^{0_f}_{c_f}\!\!&\!^{1_f}_{d_f}\!\!&\!^{2_f}_{e_f}\!\!&\!^{3_f}_{-}\!\!&\!^{4_f}_{-}\!\!&\!^{5_f}_{-}\!\!&\!^{6_f}_{-}\!\!&\!^{7_f}_{-}\!\!&\!^{8_f}_{-}\!\!&\!^{9_f}_{-}\!\!&\!^{a_f}_{-}\!\!&\!^{b_f}_{-}\\
^{11}_{12}\!\!&\!^{-}_{3_9}\!\!&\!^{-}_{4_9}\!\!&\!^{-}_{5_9}\!\!&\!^{0_9}_{6_9}\!\!&\!^{1_9}_{7_9}\!\!&\!^{2_9}_{8_9}\!\!&\!^{02}_{12}\!\!&\!^{c_f}_{-}\!\!&\!^{d_f}_{-}\!\!&\!^{e_f}_{-}\!\!&\!^{-}_{0_f}\!\!&\!^{-}_{1_f}\!\!&\!^{-}_{2_f}
\!&\!\!^{11}_{12}\!\!&\!^{-}_{9_f}\!\!&\!^{-}_{a_f}\!\!&\!^{-}_{b_f}\!\!&\!^{0_f}_{c_f}\!\!&\!^{1_f}_{d_f}\!\!&\!^{2_f}_{e_f}\!\!&\!^{3_f}_{-}\!\!&\!^{4_f}_{-}\!\!&\!^{5_f}_{-}\!\!&\!^{6_f}_{-}\!\!&\!^{7_f}_{-}\!\!&\!^{8_f}_{-}\\
\!\!&\!\!\!&\!\!\!&\!\!\!&\!\!\!&\!\!\!&\!\!\!&\!^{13}_{14}\!\!&\!^{3_f}_{9_f}\!\!&\!^{4_f}_{a_f}\!\!&\!^{5_f}_{b_f}\!\!&\!^{6_f}_{c_f}\!\!&\!^{7_f}_{d_f}\!\!&\!^{8_f}_{e_f}
\!&\!\!^{22}_{23}\!\!&\!^{-}_{6_f}\!\!&\!^{-}_{7_f}\!\!&\!^{-}_{8_f}\!\!&\!^{-}_{9_f}\!\!&\!^{-}_{a_f}\!\!&\!^{-}_{b_f}\!\!&\!^{0_f}_{c_f}\!\!&\!^{1_f}_{-}\!\!&\!^{2_f}_{-}\!\!&\!^{3_f}_{-}\!\!&\!^{4_f}_{-}\!\!&\!^{5_f}_{-}\\
\!\!\!\!&\!\!\!&\!\!\!&\!\!\!&\!\!\!&\!\!\!&\!\!\!&\!\!\!&\!\!\!&\!\!\!&\!\!\!&\!\!\!&\!\!\!&\!
\!&\!\!^{33}_{34}\!\!&\!^{-}_{3_f}\!\!&\!^{-}_{4_f}\!\!&\!^{-}_{5_f}\!\!&\!^{-}_{6_f}\!\!&\!^{-}_{7_f}\!\!&\!^{-}_{8_f}\!\!&\!^{-}_{9_f}\!\!&\!^{-}_{a_f}\!\!&\!^{-}_{b_f}\!\!&\!^{0_f}_{c_f}\!\!&\!^{1_f}_{d_f}\!\!&\!^{2_f}_{e_f}
\end{array}$$

Below the column-headers row, each such table starts  with the row-header $\alpha\beta=00$ and then cites the vertices $1_m,2_m,\ldots,(n-1)_m$ of $Y$ under the successive column headers, from left to right. Successive rows contain the remaining vertices $n_m,\ldots,(m-1)_m$\,, again from left to right. If this process does not reach the last column after setting the vertex $(m-1)_m$ of $Y$, then we set hyphens in the remaining places of the row, indicating inadequacy for further assignments there. The bicolor $\alpha'\beta'$ heading each row, except the initial row (headed by $\alpha\beta$), satisfies the following relations with respect to the bicolor $\alpha\beta$ of the preceding row: $\alpha'=\alpha$ and $\beta'=\beta+1$. If a row
has a hyphen in the last column, then its header bicolor $\alpha\beta$ has its subsequent bicolor $\alpha''\beta''$ with $\alpha''=w+1$ and $\beta''=\beta$. This process is repeated until the last column is assigned $(m-1)_m$.

\vspace*{2mm}

The algorithm allows to edit a guide for a $(\lambda,\mu)$-coloring for $H$. To obtain such a guide, we first define the $(m,n)$-{\it list} associated to $H$ as a sequence $(A_0,\ldots,A_s)$ of positive integer sequences $A_i=(a_0^i,a_1^i,\ldots,a_{r_i}^i)=A_i^{r_i}$ such that $m=a_0^i+a_1^i+\cdots+a_{r_i}^i$, for $i=0,\ldots,s$, with $n=a_0^0=a_{r_0}^0+a_0^1=\ldots=a_{r_i}^i+a_0^{i+1}=\ldots=a_{r_{s-1}}^{s-1}+a_0^s=a_{r_s}^s$
and each $a_i^j$ selected as large as possible once the set $\{a_{i'}^{j'};j'<j,\mbox{ and }i'<i\mbox{ if }j'=j\}$ is determined.

To define $(\lambda,\mu)$-guides, we need just $(\lambda,\mu)$-lists, some of which are listed in the following table:

$$\begin{array}{||l|l||l|l||}
^{\lambda,\mu}_{==}&^{A_0,A_1,\ldots,A_s}_{==============}&^{\lambda,\mu}_{==}&^{A_0,A_1,\ldots,A_s}_{=================}\\
^{4,3}_{5,2}&^{(3,1),(2,2),(1,3)}_{(2,2,1),(1,2,2)}&^{7,2}_{7,3}&^{(2,2,2,1),(1,2,2,2)}_{(3,3,1),(2,3,2),(1,3,3)}\\
^{5,3}_{5,4}&^{(3,2),(1,3,1),(2,3)}_{(4,1),(3,2),(2,3),(1,4)}&^{7,4}_{7,5}&^{(4,3),(1,4,2),(2,4,1),(3,4)}_{(5,2),(3,4),(1,5,1),(4,3),(2,5)}\\
^{6,5}&^{(5,1),(4,2),(3,3),(2,4),(1,5)}&^{7,6}&^{(6,1),(5,2),(4,3),(3,4),(2,5),(1,6)}
\end{array}$$

To each $(\lambda,\mu)$-list $(A_0^{r_0},\ldots,A_s^{r_s})$ we associate the $(\lambda,\mu)$-{\it guide}
$(B_0,\ldots,$ $B_s)$ given by replacing $A_0$ by the bicolor sequence $B_0=((0,0),\ldots,(0,r_0))$; then $A_1$ by $B_1=((1,r_0),(1,r_0+1),\ldots,(1,r_0+r_1))$; $\ldots$; then $A_i$ by $B_i=((i,r_0+\cdots+r_{i-1}),\ldots,(i,r_0+\cdots+r_{i-1}+r_i))$; $\ldots$; then $A_s$ by $B_s=((s,r_0,+\cdots+r_{s-1}),\ldots,(s,r_0+\cdots+r_{s-1}+r_s))$.
The $(\lambda,\mu)$-guides for the above $(\lambda,\mu)$-lists are as follows, where $\alpha\beta=(\alpha,\beta)$ expresses each entry of the $B_i$:

$$\begin{array}{||l|l||l|l||}
^{\lambda,\mu}_{==}&^{B_0,B_1,\ldots,B_s}_{=============
========}&^{\lambda,\mu}_{==}&^{B_0,B_1,\ldots,B_s}_{=========================}\\
^{4,3}_{5,2}&^{(00,01),(11,12),(22,23)}_{(00,01,02),(12,13,14)}&^{7,2}_{7,3}&^{(00,01,02,03),(13,14,15,16)}_{(00,01,02),(12,13,14),(24,25,26)}\\
^{5,3}_{5,4}&^{(00,01),(11,12,13),(23,33)}_{(00,01),(11,12),(22,23),(33,34)}&^{7,4}_{7,5}&^{(00,01),(11,12,13),(23,24,25),(35,36)}_{(00,01),(11,12),(22,23,24),(34,35),(45,46)}\\
^{6,5}&^{(00,01),(11,12),(22,23),(33,34),(44,45)}&^{7,6}&^{(00,01),(11,12),(22,23),(33,34),(44,45),(55,56)}
\end{array}$$

\begin{coro}\label{co1}
The $(\lambda,\mu)$-coloring in Theorem~\ref{thm1} and its associated palette $\Theta$ are in one-to-one correspondence respectively with the $(m,n)$-list of $H$ and with the $(\lambda,\mu)$-guide obtained from the corresponding $(\lambda,\mu)$-list.
\end{coro}

\proof
The $(m,n)$-list of $H$ and $ (\lambda,\mu)$-guide in the statement were conceived in relation to the construction yielded by the algorithm in the proof of Theorem~\ref{thm1}, so they are indeed in one-to-one correspondence with the $(\lambda,\mu)$-coloring and its associated palellte $\Theta$.
\qfd

\noindent The pairs $(m,n)$ of the three examples above, namely $(18,12)$, $(15,6)$ and $(15,12)$, have respective greatest common divisors $gcd(m,n)=d=6,3,3$.
The three correspond to the coprime pairs $(\lambda,\mu)=(3,2),(5,2),(5,4)$, where $\lambda=\frac{m}{d}$ and $\mu=\frac{n}{d}$.
For them,
with $d=1$, $m=\lambda$ and $n=\mu$, the following triple-table display (completed with two examples, see below) is obtained:

$$\begin{array}{||c|cc||c|cc||c|cccc||c|ccc||c|cccccc||}
^{\alpha\beta}_{==}\!&\!^{0_2}_{=}\!&\!^{1_2}_{=}\!&\!^{\alpha\beta}_{==}\!&\!^{0_2}_{=}\!&\!^{1_2}_{=}
\!&\!^{\alpha\beta}_{==}\!&\!^{0_4}_{=}\!&\!^{1_4}_{=}\!&\!^{2_4}_{=}\!&\!^{3_4}_{=}\!&\!
^{\alpha\beta}_{==}\!&\!^{0_3}_{=}\!&\!^{1_3}_{=}\!&\!^{2_3}_{=}\!&\!^{\alpha\beta}_{==}\!&\!^{0_3^0}_{=}\!&\!^{0_3^1}_{=}\!&\!^{1_3^0}_{=}\!&\!^{1_3^1}_{=}\!&\!^{2_3^0}_{=}\!&\!^{2_3^1}_{=}\\
^{00}_{01}\!&\!^{0_3}_{2_3}\!&\!^{1_3}_{-}\!&\!^{00}_{01}\!&\!^{0_5}_{2_5}\!&\!^{1_5}_{3_5}
\!&\!^{00}_{01}\!&\!^{0_5}_{4_5}\!&\!^{1_5}_{-}\!&\!^{2_5}_{-}\!&\!^{3_5}_{-}\!&\!
^{00}_{01}\!&\!^{0_5}_{3_5}\!&\!^{1_5}_{4_5}\!&\!^{2_5}_{-}\!&\!^{00}_{01}\!&\!^{0_5^0}_{3_5^0}\!&\!^{0_5^1}_{3_5^1}\!&\!^{1_5^0}_{4_5^0}\!&\!^{1_5^1}_{4_5^1}\!&\!^{2_5^0}_{-}\!&\!^{2_5^1}_{-}\\
^{11}_{12}\!&\!^{-}_{1_3}\!&\!^{0_3}_{2_3}\!&\!^{02}_{12}\!&\!^{4_5}_{-}\!&\!^{-}_{0_5}
\!&\!^{11}_{12}\!&\!^{-}_{3_5}\!&\!^{0_5}_{4_5}\!&\!^{1_5}_{-}\!&\!^{2_5}_{-}\!&\!
^{11}_{12}\!&\!^{-}_{1_5}\!&\!^{-}_{2_5}\!&\!^{0_5}_{3_5}\!&\!^{11}_{12}\!&\!^{-}_{1_5^0}\!&\!^{-}_{1_5^1}\!&\!^{-}_{2_5^0}\!&\!^{-}_{2_5^1}\!&\!^{0_5^0}_{3_5^0}\!&\!^{0_5^1}_{3_5^1}
\vspace*{0.25mm}\\
\!&\!\!&\!\!&\!^{13}_{14}\!&\!^{1_5}_{3_5}\!&\!^{2_5}_{4_5}
\!&\!^{22}_{23}\!&\!^{-}_{2_5}\!&\!^{-}_{3_5}\!&\!^{0_5}_{4_5}\!&\!^{1_5}_{-}\!&\!
^{13}_{23}\!&\!^{4_5}_{-}\!&\!^{-}_{0_5}\!&\!^{-}_{1_5}\!&\!^{13}_{23}\!&\!^{4_5^0}_{-}\!&\!^{4_5^1}_{-}\!&\!^{-}_{0_5^0}\!&\!^{-}_{0_5^1}\!&\!^{-}_{1_5^0}\!&\!^{-}_{1_5^1}
\vspace*{0.25mm}\\
\!&\!\!&\!\!&\!\!&\!\!&\!\!&\!^{33}_{34}\!&\!^{-}_{1_5}\!&\!^{-}_{2_5}\!&\!^{-}_{3_5}\!&\!^{0_5}_{4_5}\!&\!
^{24}\!&\!^{2_5}\!&\!^{3_5}\!&\!^{4_5}\!&\!^{24}\!&\!^{2_5^0}\!&\!^{2_5^1}\!&\!^{3_5^0}\!&\!^{3_5^1}\!&\!^{4_5^0}\!&\!^{4_5^1}
\end{array}$$

The $(m,n)$-lists for our three initial examples (of the algorithm), followed by the $(\lambda,\mu)$-lists of  the last five (with their $(\lambda,\mu)$-guides cited below) are:

\begin{equation}^{((6,3),(3,6)),\,((6,6,3),(3,6,6)),\,((12,3),(9,6),(6,9),(3,12)),}_{((2,1),(1,2)),\;((2,2,1),(1,2,2)),\;((4,1),(3,2),(2,3),(1,4)),\;((3,2),(1,3,1),(2,3)),\;((6,4),(2,6,2),(4,6)),}\end{equation} 
$$(^{((00,01),(11,12)),\,((00,01,02),(12,13,14)),\,((00,01),(11,12),(22,23),(33,34)),\mbox{ (these, twice each) }}_{((00,01),(11,12,13),(23,24)),\,((00,01),(11,12,13),(23,24)).}),$$

\noindent the last $(\lambda,\mu)$-list here obtained by duplicating extension of the integers of the penultimate one.
These two rightmost tables in the quintuple-table display are first for $(m,n,d)=(\lambda,\mu,1)=(5,3,1)$ and then for $(m,n,d)=(2\lambda,2\mu,2)=(10,6,2)$. The  graphs corresponding to these two tables are to be denot ed $H$ and $G$ in Section~\ref{s4} (with $H$ already cited by the end of Section~\ref{s1}).
Such a duplication extends to $d$-plication, or $d$-replication, $(1<d\in\Z$), as follows.

\begin{coro}\label{co2} Let $m=\lambda d$ and $n=\mu d$, ($1<d\in\Z$). Then, a $(\lambda,\mu)$-biregular graph $G$ with bipartition $(Y,X)$ such that $|Y|=n$ and $|X|=m$ has a reversible full $(\lambda,\mu)$-coloring. Moreover, $G$ is a $d$-ple (double, if $d=2$) cover of a $(\lambda,\mu)$-biregular graph $H$ with bipartition $(Y_H,X_H)$ with $|Y_H|=\mu$ and $|X_H|=\lambda$ that has a reversible full $(\lambda,\mu)$-coloring as in Theorem~\ref{thm1}.
\end{coro}

\proof The elements
$(di)_m,(d(i+1))_m,\ldots,(d(i+d-1))_m$ of $\Z_m$ and $(di)_n,$ $(d(i+1))_n,\ldots,(d(i+d-1))_n$ of $\Z_n$ can be taken as the vertices
$i_{\lambda}^0,i_{\lambda}^1,\ldots,$ $i_{\lambda}^{d-1}$ and $i_{\mu}^0,i_{\mu}^1,\ldots,$ $i_{\mu}^{d-1}$ of $G$, respectively, obtained by $d$-plication from 
the vertices $i_{\lambda}\in\Z_{\lambda}$ and $i_{\mu}\in\Z_{\mu}$ of $H$.
These elements are easily seen to have a disposition as above, via a straightforward generalization of that of the rightmost table in the quintuple-table display (for $G$) with respect to the one at its immediate left (for $H$). 
There, instead of the elements $(2i)_{10},(2i+1)_{10}$, ($i\in \Z_5$)  and $(2i)_6,(2i+1)_6$, ($i\in \Z_3$) of our previous notation in $H$, we use $i_5^0,i_5^1$ and $i_3^0,i_3^1$, respectively. This table yields a reversible full $(5,3)$-coloring of $G$. Similarly in the general situation of the statement of the corollary.
\qfd

The weights of the pairs produced by the algorithm constitute in their order of appearance the {\it weight distribution} of the $(\lambda,\mu)$-coloring.
For example, the weight distributions in the first display above and then in the recent display are:

$$^{(6,3,3,6),(6,6,3,3,6,6),(12,9,3,6,6,3,9,12);}_{(2,1,1,2),(2,2,1,1,2,2),(4,\hspace*{0.5mm}3,\hspace*{0.5mm}1,2,\hspace*{0.5mm}2,\hspace*{0.5mm}1,\hspace*{0.5mm}3,\hspace*{0.5mm}4), (3,2,1,3,1,2,3),(6,4,2,6,2,4,6).}$$

\section{Odd-Girth Cycles in Odd Graphs}\label{s5/2}

Let $K_n$ be the complete graph with vertex set $[n]=\{1,2,\ldots,n\}$, ($3\le n\in\Z$). First, we establish some notation for $n$-cycles (of $K_n$).
Every permutation of $[n]$ can be written as an ordered list $\ell_n=a_1\cdots a_n$ of pairwise distinct integers $a_i\in[n]$, ($i\in[n]$). All such lists $\ell_n$ are distributed into $(n-1)!$ classes, where each class: {\bf(a)} represents $n$ distinct permutations of $[n]$ and {\bf(b)} is representable as an $n$-circuit $(a_1,\ldots,a_n)$ of $K_n$, 
i.e. a subdigraph of $K_n$ with arcs $(a_1,a_2),(a_2,a_3),\ldots,$ $(a_{n-1},a_n),$ $(a_n,a_1)$. By disregarding the orientations of these $n$-circuits, we are left with $\frac{1}{2}(n-1)!$ $n$-cycles of $K_n$ that we denote $(\ell_n)= (a_1\ldots a_n)$, between parentheses and with no commas, the cited arcs replaced by the corresponding edges $a_1a_2,a_2a_3,\ldots,$ $a_{n-1}a_n,a_na_1$.

From now on, let $n=2k+1>3,k\in\Z$. The odd graph $O_{k+1}$ \cite{B,BCN,GR} can be defined to have the $(k+1)$-subsets of $[n]$ as its vertices, and an edge $e$ between two such vertices if and only if
the $(k+1)$-subsets they represent intersect at just one element $i\in[n]$, so $e$ is assigned color $i$. This way, $O_{k+1}$ becomes edge colored. For example, each vertex of the Petersen graph $O_3$ is denoted by the triple of colors assigned to its incident edges.

The following statement seems related to Corollary 7.8.2 \cite{GR}, which yields an isomorphism between the
automorphisms groups of $K_n$ and $O_{k+1}$.

\begin{thm}\label{t3} There exists a bijection $\Phi$
from the family of (Hamilton) $n$-cycles of $K_n$ onto the family of $n$-cycles of $O_{k+1}$ (which realize its odd girth).
\end{thm} 

\proof 
Each Hamilton cycle  of $K_n$ is an $n$-cycle $C$ that can be expressed $C=(a_0b_1a_1b_2a_2b_3\ldots b_ia_i\ldots b_ka_k)$, where the set $\{a_0$, $b_1$, $a_1$, $b_2$, $a_2$, $b_3$, $\ldots$, $b_i$, $a_i$, $\ldots$, $b_k$, $a_k\}$ is in one-to-one correspondence with $[n]$. For each $n$-cycle $C$ of $K_n$ written this way there exists a unique $n$-cycle $\Phi(C)$ in $O_{k+1}$ whose induced edge-color cycle $C'$ coincides with the assumed expression of $C$ above, namely:

\begin{equation}\begin{array}{ll}
^{\Phi(C)=}_{a_0b_k\cdots b_2b_1,}&^{\!\!\!(a_k\cdots a_1a_0,}_{b_1a_k\cdots a_2a_1,}\\
^{a_1a_0b_k\cdots b_3b_2,}_{a_2a_1a_0b_k\cdots b_4b_3,}&^{b_2b_1a_k\cdots a_3a_2,}_{b_3b_2b_1a_k\cdots a_4a_3,}\\
^{\cdots\cdots\cdots\cdots,}_{a_ia_{i-1}\cdots a_0b_k\cdots b_{i+2}b_{i+1},}&^{\cdots\cdots\cdots\cdots,}_{b_{i+1}b_i\cdots b_1a_k\cdots a_{i+2}a_{i+1},}\\
^{\cdots\cdots\cdots\cdots,}_{a_{k-1}\cdots a_1a_0b_k,}&^{\cdots\cdots\cdots\cdots,}_{b_kb_{k-1}\cdots b_1a_k),}
\end{array}\end{equation}

\noindent where the $n$ entries of the cycle $\Phi(C)$ in $O_{k+1}$ represents $(k+1)$-subsets of $[n]$ which are expressed without braces or parentheses and as convenient lists whose initial (or terminal) elements constitute  $C$ in the order in which they appear from left to right. (Observe that these lists are the successive $(k+1)$-sublists of the circular $n$-list $(a_k\cdots a_1a_0b_k\cdots b_2b_1)$). In fact,  
the initial element $a_k$ of $[n]$ in the expression of $\Phi(C)$ in (7) followed by the elements $a_0,b_1,\ldots,b_{k}$ of $[n]$ that follow the subsequent commas in $\Phi(C)$ form $C'$. 
Now, observe that $C$ and $C'$ have exactly the same form. Thus, $C=(a_0b_1a_1b_2a_2b_3\ldots b_ia_i\ldots b_ka_k)=C'$. This insures that $\Phi$ is a bijection.
\qfd

Note that the $(n-1)!/2$ $n$-cycles of $K_n$ are divided into two sets of $n$-cycles, namely: $(n-1)!/4$ that we call {\it evenly} $n$-cycles and are obtained from the identity-permutation $n$-cycle $(12\cdots n)$ by applying an even number of transpositions, and the remaining $(n-1)!/4$, that we call {\it oddly} $n$-cycles.

\section{Approach via the Petersen Graph}\label{s3}

From now on, let $k=2$ and $n=5$.
Let $H'$ be the $(n,3)$-biregular graph with vertex parts $Y=\{$oddly $n$ cycles of $K_n\}$ and $X=\{$3-subsets of $[n]\}$, with an edge $yx$ whenever the $n$-cycle $y$ has two contiguous entries $p,q\in x$ that are both nonadjacent in $y$ to the sole $r\in x\setminus\{p,q\}$. We properly color each edge $yx$ of $H'$ with its corresponding $r$.
 There are 6 oddly and 6 evenly 5-cycles in $O_3$, and we denote them $y_i$ and $\bar{y_i}$ ($i\in\Z_6$), respectively, with:

$$
^{y_0=(13524),\; y_1=(12435),\; y_2=(12354),\; y_3=(15234),\; y_4=(13452),\; y_5=(13245),}
_{\bar{y}_0=(12345),\;  \bar{y}_1=(13254),\; \bar{y}_2=(15243),\; \bar{y}_3=(13542),\; \bar{y}_4=(15324), \; \bar{y}_5=(14352).}$$

On the other hand, the 3-sets of $[n]=[5]$ can be considered as 3-cycles of $K_5$.
Accordingly, we denote the vertices of $O_3$ as follows:

$$^{ x_1=(314),\;\, x_2=(425),\;\, x_ 3=(531),\;\, x_4=(142),\;\, x_5=(253),}
  _{ z_1=(512),\;\; z_2=(123),\;\; z_3=(234),\;\; z_4=(345),\;\; z_5=(451).}$$

\begin{figure}[htp]\hspace{25mm}
\includegraphics[scale=0.25]{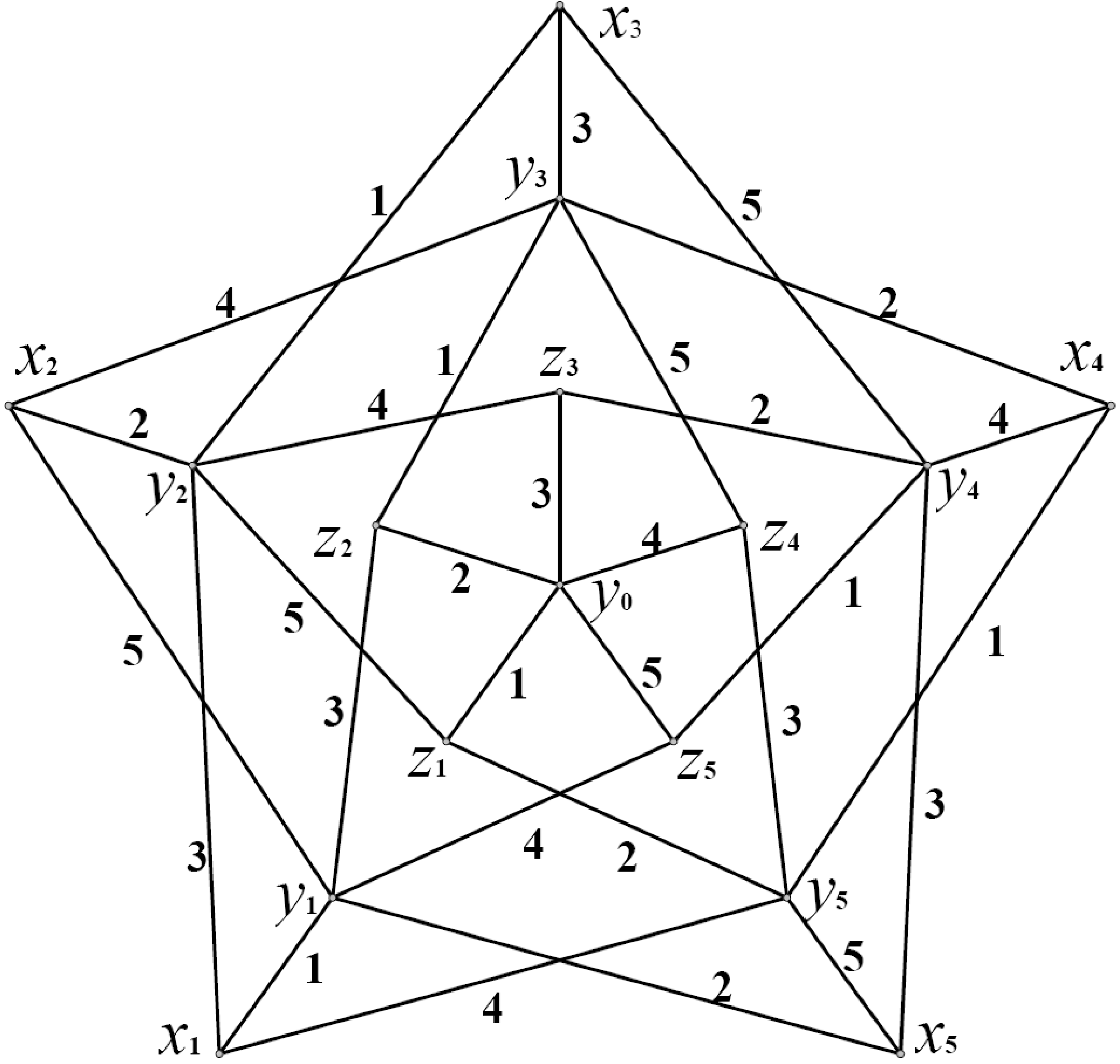}
\caption{Representation of $H'$ based on $O_3$}
\end{figure}

To focus on our assembling of the Great Circle Challenge Puzzle, $\Z_5\times\Z_3$ is redenoted $[5]\times\{a,b,c\}$ (by sending, say, $i_5\in\Z_5$ onto $i+1\in[5]$, and
$0_3,1_3,2_3\in\Z_3$ respectively onto $a,b,c$), so now $\Theta=\{1a,2a,2b,3b,4b,4c,5c\}$. 

\begin{thm}\label{abc} The graph $H'$ is a $(5,3)$-biregular graph having at least 3 full $(5,3)$-colorings with associated palette $\Theta=\{1a,2a,2b,3b,4b,4c,5c\}$. More\-over, there are at least one reversible and two non-reversible such colorings.
\end{thm}

\proof The three solutions cited at the end of Section~\ref{s1} arise as follows:

{\bf(A)} a reversible full $(5,3)$-coloring of $H'$ is obtained via colors $\alpha$ and $\beta$ associated respectively to the
degree-5 and -3 endpoints $v_e^5$ and $v_e^3$ of the edges $e$ of $H'$, and satisfying that bicolors $\alpha\beta=r\beta=$

{\bf 1.} $1a,3b,5c$ are assigned each to 6 edges of
$H'$,

{\bf 2.} $2a,4c$ are assigned each to 4 edges of
$H'$ and

{\bf 3.} $2b,4b$ are assigned each to 2 edges of $H'$.

\noindent In fact, bicolors $1a,3b,5c$ are assigned to the edges of $H'$ that in Figure 1 have respective colors 1, 3, 5; the other edges of $H'$ are assigned bicolors $\alpha\beta=r\beta$, as in the leftmost table in the triple-table display below. We note this $(5,3)$-coloring of $H'$ is reversible via the permutation $(1,5)(2,4)(a,c)$: 

{\bf(B)} A non-reversible $(5,3)$-coloring of $H'$ is obtained by assigning bicolor $\alpha\beta=r\beta=1a$ to the edges of $H'$ colored with 1, and the remaining bicolors $\alpha\beta$ as shown on the middle table of the triple-table display below, where $\beta\ne r$ happens 15 times, while in (C) below, $\beta\ne r$ happens 14 times. 

{\bf(C)} We present a third $(5,3)$-coloring of $H'$ to show that the starting and ending levels of a $(5,3)$-coloring of $H'$ can be
placed with the assignment of $1a$ and $5c$ to the edges of $H'$ colored respectively with 1 and 5, as in the
reversible $(5,3)$-coloring above, but providing the rest of the (5,3)-coloring differently,  shown on the rightmost table of the triple-table display below.\qfd

$$\begin{array}{||l|cccccc||l|cccccc||l|cccccc||}
^{(A)}_{\alpha\beta}\!&\!_{y_0}\!&\!_{y_1}\!&\!_{y_2}\!&\!_{y_3}\!&\!_{y_4}\!&\!_{y_5}\!&\!^{(B)}_{\alpha\beta}\!&\!_{y_0}\!&\!_{y_1}\!&\!_{y_2}\!&\!_{y_3}\!&\!_{y_4}\!&\!_{y_5}\!&\!^{(C)}_{\alpha\beta}\!&\!_{y_0}\!&\!_{y_1}\!&\!_{y_2}\!&\!_{y_3}\!&\!_{y_4}\!&\!_{y_5}\\
^{==}_{1a}\!&\!^{=}_{z_1}\!&\!^{=}_{x_1}\!&\!^{=}_{x_3}\!&\!^{=}_{z_2}\!&\!^{=}_{z_5}\!&\!^{=}_{x_4}\!&\!^{==}_{1a}\!&\!^{=}_{z_1}\!&\!^{=}_{x_1}\!&\!^{=}_{x_3}\!&\!^{=}_{z_2}\!&\!^{=}_{z_5}\!&\!^{=}_{x_4}\!&\!^{==}_{1a}\!&\!^{=}_{z_1}\!&\!^{=}_{x_1}\!&\!^{=}_{x_3}\!&\!^{=}_{z_2}\!&\!^{=}_{z_5}\!&\!^{=}_{x_4}\\
^{2a}_{2b}\!&\!^{z_2}_{-}\!&\!^{x_5}_{-}\!&\!^{x_2}_{-}\!&\!^{x_4}_{-}\!&\!^{-}_{z_3}\!&\!^{-}_{z_1}\!&\!^{2a}_{2b}\!&\!^{z_3}_{-}\!&\!^{-}_{x_5}\!&\!^{x_2}_{-}\!&\!^{-}_{x_4}\!&\!^{x_5}_{-}\!&\!^{z_4}_{-}\!&\!^{2a}_{2b}\!&\!^{-}_{z_3}\!&\!^{x_5}_{-}\!&\!^{x_2}_{-}\!&\!^{-}_{x_3}\!&\!^{z_3}_{-}\!&\!^{z_4}_{-}\\
^{3b}_{4b}\!&\!^{z_3}_{z_4}\!&\!^{z_2}_{z_5}\!&\!^{x_1}_{-}\!&\!^{x_3}_{-}\!&\!^{x_5}_{-}\!&\!^{z_4}_{-}\!&\!^{3b}_{4b}\!&\!^{z_4}_{-}\!&\!^{z_2}_{z_5}\!&\!^{x_1}_{-}\!&\!^{x_2}_{-}\!&\!^{z_3}_{x_3}\!&\!^{z_1}_{-}\!&\!^{3b}_{4b}\!&\!^{z_2}_{z_4}\!&\!^{z_5}_{-}\!&\!^{x_1}_{-}\!&\!^{x_4}_{x_2}\!&\!_{-}^{x_5}\!&\!_{-}^{z_1}\\
^{4c}_{5c}\!&\!^{-}_{z_5}\!&\!^{-}_{x_2}\!&\!^{z_3}_{z_1}\!&\!^{x_2}_{z_4}\!&\!^{x_4}_{x_3}\!&\!^{x_1}_{x_5}\!&\!^{4c}_{5c}\!&\!^{z_2}_{z_5}\!&\!^{-}_{x_2}\!&\!^{z_1}_{z_3}\!&\!^{x_3}_{z_4}\!&\!^{-}_{x_4}\!&\!^{x_5}_{x_1}\!&\!^{4c}_{5c}\!&\!^{-}_{z_5}\!&\!^{z_2}_{x_2}\!&\!^{z_3}_{z_1}\!&\!^{-}_{z_4}\!&\!^{x_4}_{x_3}\!&\!^{x_1}_{x_5}\\
\end{array}$$

\section{Duplicating Extension}\label{s4}

Recall the reversible full $(5,3)$-coloring of the $(5,3)$-bi\-re\-gu\-lar graph $H'$  in Theorem~\ref{abc}. We will build now the  bipartite double cover $G'$ of $H'$.
We note that $H'$ has $(m,n,d)=(10,6,2)$ while $G'$ will be shown to have $(m,n,d)=(20,12,4)$. The weight distribution of $H'$ is $(6,2,4,6,4,2,6)$ while that of $G'$ will be  $(12,4,8,12,8,4,12)$.

\begin{thm}\label{assem}
There exists a bipartite double cover $G'$ of $H'$ that is a $(5,3)$-biregular graph. Moreover, $G'$ has at least three full $(5,3)$-colorings with associated palette $\Theta=\{1a,2a,2b,3b,4b,4c,5c\}$. Furthermore, there exist at least one reversible and two non-reversible such colorings.
\end{thm}

\proof We consider the face-colored plane graph in Figure 2 (with the colors of faces given in the second paragraph of Section~\ref{s5}). In the figure, a 3-path is said to be a {\it special 3-path} (or {\it s3p}) if it is drawn with:

{\bf(a)} a thick middle edge;

{\bf(b)} two dashed end edges.

\noindent Let $\cal G$ be the graph obtained from the plane graph in Figure 2 by 
identifying the 4 corner vertices labeled $1_5$
as a result of folding along the 4 hypotenuses $1_41_3$, $1_61_4$, $1_11_6$, $1_31_1$ of the isosceles right triangles:

\begin{figure}[htp]
\hspace*{16mm}
\includegraphics[scale=0.59]{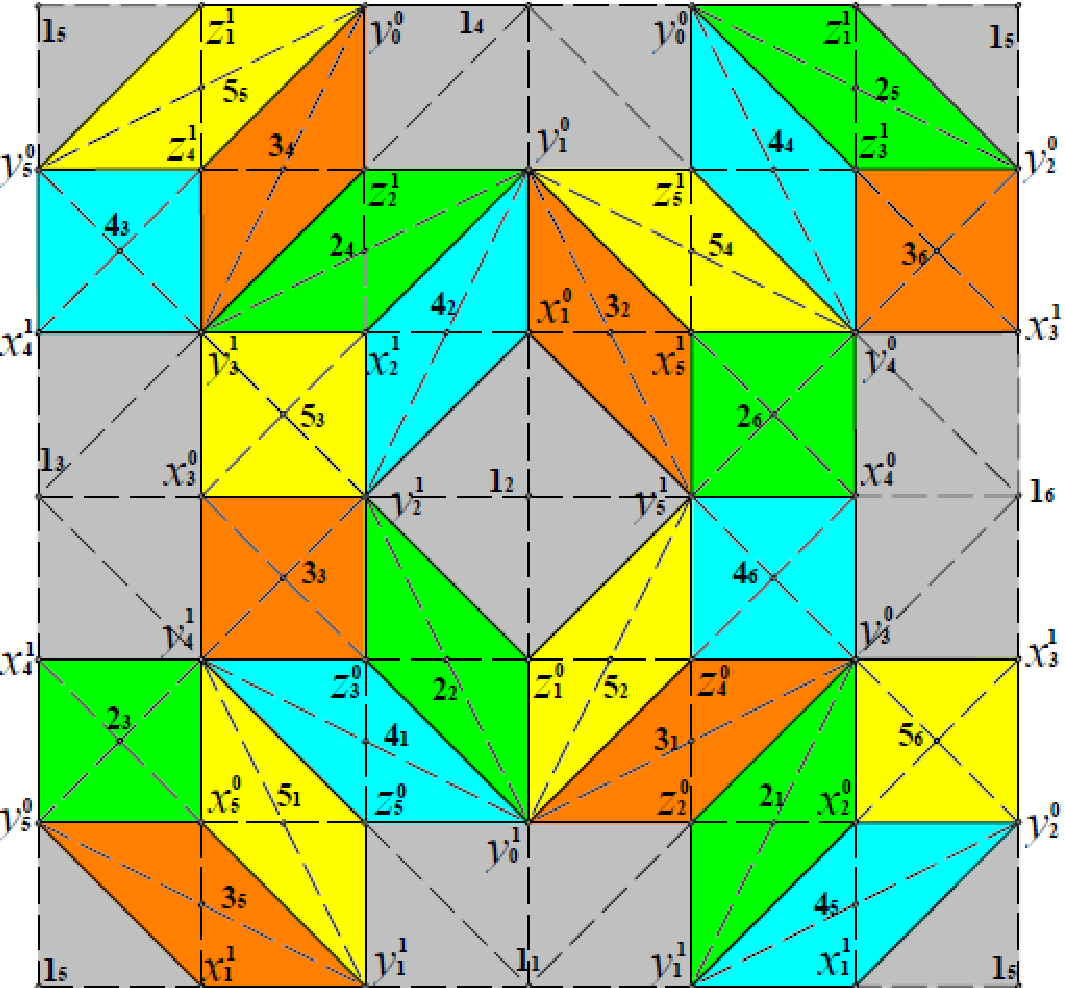}
\caption{Cutout of the face-colored subgraph $G'$ of graph $\cal G$}
\end{figure}

$$\begin{array}{cc}
^{\mbox{\tiny{triangle}}(1_51_41_3)=\{1_5z_1^1y_0^01_4,\;1_4z_2^1y_3^11_3,\;1_3x_4^1y_5^01_5\},}
_{\mbox{\tiny{triangle}}(1_21_41_3)=\{1_2 x_1^0 y_1^0 1_4 ,\;1_4 z_2^1 y_3^1 1_3 ,\;1_3 x_3^0 y_4^1 1_2 \},}\!&\!
^{\mbox{\tiny{triangle}}(1_51_61_4)=\{1_5y_2^0x_3^11_6,\;1_6y_4^0z_5^11_4,\;1_4y_0^0z_1^11_5\},}
_{\mbox{\tiny{triangle}}(1_21_61_4)=\{1_2 y_5^1 x_4^0 1_6 ,\;1_6 y_4^0 z_5^1 1_4 ,\; 1_4 y_1^0 x_1^0 1_2 \},}\\
^{\mbox{\tiny{triangle}}(1_51_11_6)=\{1_5x_1^1y_1^11_1,\;1_1z_2^0y_3^01_6,\;1_6x_3^1y_2^01_5\},}
_{\mbox{\tiny{triangle}}(1_21_11_6)=\{1_2 z_1^0 y_0^1 1_1 ,\; 1_1 z_2^0 y_3^0 1_6 ,\; 1_6 x_4^0 y_5^1 1_2 \},}\!&\!
^{\mbox{\tiny{triangle}}(1_51_31_1)=\{1_5y_5^0x_4^11_3,\;1_3y_4^1z_5^01_1,\;1_1y_1^1x_1^11_5\},}
_{\mbox{\tiny{triangle}}(1_21_31_1)=\{1_2 y_2^1 x_3^0 1_3, \; 1_3 y_4^1 z_5^0 1_1 ,\;1_1 y_0^1 z_1^0 1_2 \},}
\end{array}$$

\noindent towards the central vertex (labeled) $1_2$ (but keeping $1_5\ne 1_2$ as vertices in $\cal G$) as well as identifying all other pairs of vertices with a common label 
(specifically $z_1^1,y_0^0,y_4^0,x_3^1,x_1^1,y_3^1,y_2^0,x_4^1$)
and pairs of edges joining them. This can be done simultaneously while ``inflating'', that is modifying the said 8 triangles 
into 8 corresponding spherical equilateral right triangles composing a sphere $S$ (that reappears in Section~\ref{s5} so that the external skeleton resulting from assembling the GCCP fits tightly into it). As a result, all vertex labels of $\cal G$ are different, yielding a vertex notation in $\cal G$.  

Each edge as in (b) above joins an endvertex $w_i$ of the s3p $P$ that contains it, where $w\in\{x,z\}\in[5]$ and $i\in [6]$, and a vertex of a monochromatic face (of a thick-edge 4-cycle) containing $w_i$ in its interior. Here symbol $w$ is shared by both endvertices of $P$.

Let $G'$ be the thick-edge subgraph of $\cal G$. Then, $H'$ is a quotient graph of $G'$ under a $\Z_2$-action on $\cal G$ that identifies vertex pairs (including those in ${\cal G}-G'$) of the form $\{i_j,i_{j+3}\}$, ($i\in [5]$, $j\in [3]$), and the s3ps between them.
Moreover, there is a covering graph map $\phi:G'\rightarrow H'$, with the inverse images of the vertices $y_i$ $(i\in\Z_6)$ and $x_i,z_i$ $(i\in [5])$ of $H'$
formed by the pairs of vertices $\{ y_i^0, y_i^1\}$ $(i\in\Z_6)$ and $\{ x_i^0, x_i^1\}$ $\{ z_i^0, z_i^1\}$
$(i\in [5])$ given as follows:

$$\begin{array}{llll}
^{y_0^0=(1_4,\;3_4,\;5_5,\;2_5,\;4_4),}_{y_1^0=(1_4,\;5_4,\;3_2,\;4_2,\;2_4),}&
^{y_0^1=(1_1,\;4_1,\;2_2,\;5_2,\;3_1);}_{y_1^1=(1_1,\;2_1,\;4_5,\;3_5,\;5_1);}&
^{y_3^0=(1_6,\;5_6,\;2_1,\;3_1,\;4_6),}_{y_4^0=(1_3,\;3_3,\;4_1,\;5_1,\;2_3),}&
^{y_3^1=(1_3,\;4_3,\;3_4,\;2_4,\;5_3);}_{y_4^1=(1_6,\;2_6,\;5_4,\;4_4,\;3_6);}\\
^{y_2^0=(1_5,\;4_5,\;5_6,\;3_6,\;2_5),}&
^{y_2^1=(1_2,\;2_2,\;3_3,\;5_3,\;4_2);}&
^{y_5^0=(1_5,\;5_5,\;4_3,\;2_3,\;3_5),}&
^{y_5^1=(1_2,\;3_2,\;2_6,\;4_6,\;5_2);}\\
^{x_1^0=(1_2,\;4_2,\;3_2),}_{x_2^0=(2_1,\;5_6,\;4_5),}&
^{x_1^1=(1_5,\;3_5,\;4_5);}_{x_2^1=(2_4,\;4_2,\;5_3);}&
^{z_1^0=(1_2,\;5_2,\;2_2),}_{z_2^0=(1_1,\;3_1,\;2_1),}&
^{z_1^1=(1_5,\;2_5,\;5_5);}_{z_2^1=(1_4,\;2_4,\;3_4);}
\vspace*{1mm}\\
^{x_3^0=(1_3,\;5_3,\;3_3),}_{x_4^0=(1_6,\;4_6,\;2_6),}&
^{x_3^1=(1_6,\;3_6,\;5_6);}_{x_4^1=(1_3,\;2_3,\;4_3);}&
^{z_3^0=(2_2,\;4_1,\;3_3),}_{z_4^0=(3_1,\;5_2,\;4_6),}&
^{z_3^1=(2_5,\;3_6,\;4_4);}_{z_4^1=(3_4,\;4_3,\;5_5);}
\vspace*{1mm}\\
^{x_5^0=(2_3,\;5_1,\;3_5),}&
^{x_5^1=(2_6,\;3_2,\;5_4);}&
^{z_5^0=(1_1,\;5_1,\;4_1),}&
^{z_5^1=(1_4,\;4_4,\;5_4).}
\end{array}$$

We reproduce now the leftmost table at the triple-table display for solution (A) in Section~\ref{s3}, and present to its right the one obtainable from it by duplicating extension in the bipartite double cover $G'$, yielding a reversible full $(5,3)$-coloring in $G'$:

$$\begin{array}{||c|cccccc||c||cccccccccccc||}
^{\alpha\beta}_{==}\!&\!^{y_0}_{=}\!&\!^{y_1}_{=}\!&\!^{y_2}_{=}\!&\!^{y_3}_{=}\!&\!^{y_4}_{=}\!&\!^{y_5}_{=}\!&\!
^{\alpha\beta}_{==}\!&\!
^{y_0^0}_{=}\!&\!^{y_0^1}_{=}\!&\!
^{y_1^0}_{=}\!&\!^{y_1^1}_{=}\!&\!
^{y_2^0}_{=}\!&\!^{y_2^1}_{=}\!&\!
^{y_3^0}_{=}\!&\!^{y_3^1}_{=}\!&\!
^{y_4^0}_{=}\!&\!^{y_4^1}_{=}\!&\!
^{y_5^0}_{=}\!&\!^{y_5^1}_{=}
\vspace*{0.5mm}\\
^{1a}_{2a}\!&\!
^{z_1}_{z_2}\!&\!
^{x_1}_{x_5}\!&\!
^{x_3}_{x_2}\!&\!
^{z_2}_{x_4}\!&\!
^{z_5}_{-}\!&\!
^{x_4}_{-}\!&\!
^{1a}_{2a}\!&\!
^{z_1^1}_{z_2^1}\!&\!\
^{z_1^0}_{z_2^0}\!&\!
^{x_1^0}_{x_5^1}\!&\!
^{x_1^1}_{x_5^0}\!&\!
^{x_3^1}_{x_2^0}\!&\!
^{x_3^0}_{x_2^1}\!&\!
^{z_2^0}_{x_4^0}\!&\!
^{z_2^1}_{x_4^1}\!&\!
^{z_5^1}_{-}\!&\!
^{z_5^0}_{-}\!&\!
^{x_4^1}_{-}\!&\!
^{x_4^0}_{-}
\vspace*{0.5mm}\\
^{2b}_{3b}\!&\!
^{-}_{z_3}\!&\!
^{-}_{z_2}\!&\!
^{-}_{x_1}\!&\!
^{-}_{x_3}\!&\!
^{z_3}_{x_5}\!&\!
^{z_1}_{z_4}\!&\!
^{2b}_{3b}\!&\!
^{-}_{z_3^1}\!&\!
^{-}_{z_3^0}\!&\!
^{-}_{z_2^1}\!&\!
^{-}_{z_2^0}\!&\!
^{-}_{x_1^1}\!&\!
^{-}_{x_1^0}\!&\!
^{-}_{x_3^1}\!&\!
^{-}_{x_3^0}\!&\!
^{z_3^1}_{x_5^1}\!&\!
^{z_3^0}_{x_5^0}\!&\!
^{z_1^1}_{z_4^1}&
^{z_1^0}_{z_4^0}
\vspace*{0.5mm}\\
^{4b}_{4c}\!&\!
^{z_4}_{-}\!&\!
^{z_5}_{-}\!&\!
^{-}_{z_3}\!&\!
^{-}_{x_2}\!&\!
^{-}_{x_4}\!&\!
^{-}_{x_1}\!&\!
^{4b}_{4c}\!&\!
^{z_4^1}_{-}\!&\!
^{z_4^0}_{-}\!&\!
^{z_5^1}_{-}\!&\!
^{z_5^0}_{-}\!&\!
^{-}_{z_3^1}\!&\!
^{-}_{z_3^0}\!&\!
^{-}_{x_2^0}\!&\!
^{-}_{x_2^1}\!&\!
^{-}_{x_4^0}\!&\!
^{-}_{x_4^1}\!&\!
^{-}_{x_1^1}\!&\!
^{-}_{x_1^0}
\vspace*{0.5mm}\\
^{5c}\!&\!^{z_5}\!&\!^{x_2}\!&\!^{z_1}\!&\!^{z_4}\!&\!^{x_3}\!&\!^{x_5}\!&\!
^{5c}\!&\!
^{z_5^1}\!&\!^{z_5^0}\!&\!
^{x_2^1}\!&\!^{x_2^0}\!&\!
^{z_1^1}\!&\!^{z_1^0}\!&\!
^{z_4^0}\!&\!^{z_4^1}\!&\!
^{x_3^1}\!&\!^{x_3^0}\!&\!
^{x_5^0}\!&\!^{x_5^1}
\end{array}$$

Non-reversible duplicating for the two remaining full $(5,3)$-colorings of $G'$ at the end of Section~\ref{s3}, for solutions (B)-(C), can be verified similarly.\qfd

\section{The Great Circle Challenge Puzzle}\label{s5}

The $(\lambda,\mu)$-colorings above can be applied to the GCCP assemblage problem. This consists in assembling 60 plastic {\it GCCP-pieces} into
a structure whose external boundary represents a partition of the sphere $S$ into 120 fundamental domains of the full icosahedral group \cite{Klein}. Each of these domains is represented by a corresponding internally hollow spherical triangle (or {\it GCCP-triangle}). 
Each possible assemblage of the GCCP is bounded tightly by 15 great circles (or {\it GC}s) in $S$ determining those GCCP-triangles. Each such GC is given by the external boundary of an annulus obtained by the attachment of 4 GCCP-pieces.

The thick-edge 4-cycles of $\cal G$ have the interiors of the faces they define in Figure 2 colored according to the subindex $i$ of their central vertices $w_i$. 
So, there are 6 thick-edge 4-cycles for each one of the colors 1 (gray), 2 (green), 3 (orange), 4 (blue) and 5 (yellow). They yield a total of 30 thick-edge 4-cycles, comprising 120 triangles in $\cal G$ that correspond to the 120 domains above.

There are 12 edges in $G'$ joining two thick-edge 4-cycles with a common color $i\in[5]$, (from a total of 60 edges). They inherit color $i$ via the covering graph map $\phi:G'\rightarrow H'$. This and assigning the color $i\in[5]$ of each face in $\cal G$ to the 4 dashed edges in its interior yields a coloring of $\cal G$ in which each edge class corresponds to three mutually orthogonal GCs in $S$. 

\begin{figure}[htp]
\hspace*{1.5cm}
\includegraphics[scale=0.40]{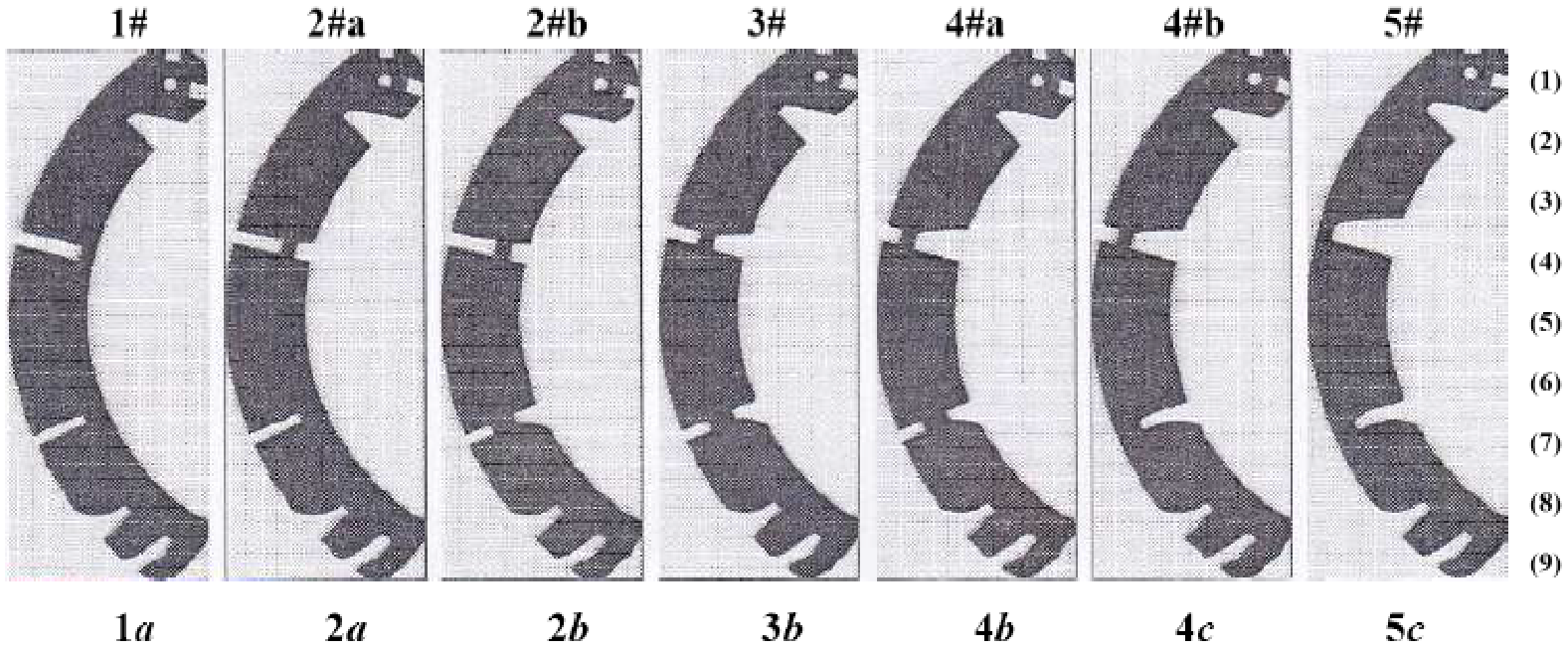}
\caption{Pieces $1^\sharp,2^\sharp$a$,2^\sharp$b$,3^\sharp,4^\sharp$a$,4^\sharp$b$,5^\sharp$, (for us: $1a,2a,2b,3b,4b,4c,5c$)}
\end{figure}

The 7 cases depicted  in Figure 3 show the front of the pairwise distinct GCCP-pieces, each covering a quarter of an annulus and spanning an angle slightly larger than $90^\circ$, measured from the center of $S$. Each GCCP-piece contains (indicated on the right of Figure 3):

{\bf(1)} a locking hole fitting a peg as in item (3) in any other GCCP-piece;

{\bf(2)} an inward cut to oppose one as in (7) for the $90^\circ$ locking of two GCs;

{\bf(3)} a short cylindrical front-side locking peg to fit as in item (1);

{\bf(4)} one or two (opposite) cuts for the successive $36^\circ$ locking of 5 GCs;

{\bf(5)} a front-side relief reading just one of $1^\sharp,2^\sharp$a$,2^\sharp$b$,3^\sharp,4^\sharp$a$,4^\sharp$b$,5^\sharp$;

{\bf(6)} one or two (opposite) cuts for the successive $60^\circ$ locking of three GCs;

{\bf(7)} a short cylindrical rotation peg on the back side;

{\bf(8)} an outward cut for the $90^\circ$ locking of two GCs, one as in (2) above;

{\bf(9)} an outward rotation cut to fit and rotate a peg as in item (7).

\noindent Pegs and cuts as in items (7) and (9) combine into hinges that allow to rotate the GCCP-pieces into correct position at their rotation
ends. Pegs and holes as in items (1) and (3) allow to attach two
GCCP-pieces at their locking ends. Two GCCP-pieces rotated and/or locked in these two ways 
show the front side of one of them attached to the back side of the
other one.

Splitting the 60 GCCP-pieces into 7 subsets distinguished as in item (5) above is controlled as follows: cases $1^\sharp$, $3^\sharp$ and $5^\sharp$ in Figure 3 (to be redenoted $1a$, $3b$ and $5c$, respectively) are given by 12 pieces each, and cases $2^\sharp$a, $2^\sharp$b, $4^\sharp$a and
$4^\sharp$b (redenoted $2a$, $2b$, $4b$ and $4c$) are given by 8, 4, 4 and 8 pieces, respectively.
These indications, $1a,2a,2b,3b,4b,4c,5c$, combine 5 cases of item (4) and three cases of item (6) above. These 8 cases are represented as cut sections of GCCP-pieces in Figure 4.

\begin{coro}\label{coro} The three full $(5,3)$-colorings of $G'$ in Theorem~\ref{assem}, that are inverse images via $\phi^{-1}$ of those corresponding to the solutions (A), (B), (C) in the proof  of Theorem~\ref{abc}, yield three corresponding solutions
to the GCCP assemblage problem, only solution (A) being reversible among them, in contrast to the always reversible situations of Theorem~\ref{thm1} and Corollaries~\ref{co1}-\ref{co2} via the cyclic-group pairs approach of Section~\ref{s2}.\end{coro}

\proof
Any 12 GCCP-pieces can be assembled into a partial structure of three annuli with external boundaries appearing
as three mutually orthogonal GCs. However, the three full $(5,3)$-colorings of Theorem~\ref{assem} of the bipartite double cover in Section~\ref{s4}  start their assemblage with the 12 GCCP-pieces $1a$.

Each GCCP-triangle has as vertices the locking intersection of: {\bf(I)} five GCs forming ten $36^\circ$ angles; {\bf(II)} three GCs forming six $60^\circ$ angles; {\bf(III)} two GCs forming four $90^\circ$ angles. There are respectively 12, 20 and 30 such vertices in
any GCCP assemblage. How do these intersections interlock?

\begin{figure}[htp]
\includegraphics[scale=0.328]{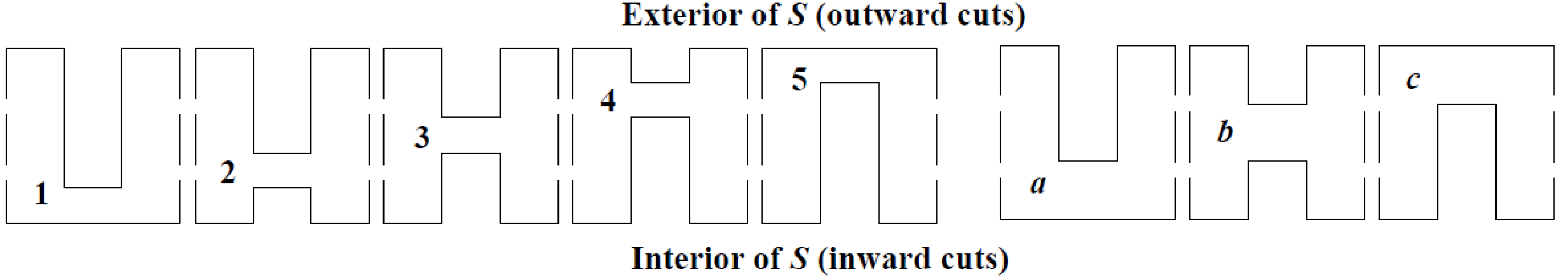}
\caption{Cut section types}
\end{figure}

To answer this, notice that the cuts in items (4) and (6) in Figure 3, provided in order to
assemble the GCCP-pieces at distinct levels from the inside out, are of two different {\it cut section types} associated to items (I) and (II) above, respectively, namely: (I) the cuts in item (4) are distributed in 5 types depicted on the left of Figure 4 and denoted $1,2,3,4,5$; (II) the cuts in item (6) are distributed in three types depicted on the right of Figure 4 and denoted $a,b,c$. The $90^\circ$ intersections as in (III) are obtained by perpendicularly setting the result of locking-and-rotating the pegs and holes about the top and bottom end parts of the GCCP-pieces.

We assign to each of these three intersection types (I), (II), (III) a
name ${\cal I}(p)$ ($p=5,3,2$, respectively) representing an intersection of $p$ GCCP-pieces
whose external borders extend to $p$ corresponding GCs.
This determines $2p$ angles of $(360/2p)^\circ$ at each intersection
point and its antipodal point in $S$.

In sum, there are: (I) 12 intersections of type ${\cal
I}(5)$; (II) 20 intersections of type ${\cal I}(3)$; and
(III) 30 intersections of type ${\cal I}(2)$. These numbers
reflect the numbers of faces, edges and vertices of the regular
dodecahedron.

Symbol 1, resp. $a$, represents a cut going deepest in $S$, so that it must be assembled first; symbol 2, resp. $b$, represents a cut one-fifth, resp. one-third, less in
depth that the one represented by symbol 1, resp. $a$; and so on
until symbol 5, resp. $c$, represents a cut to be assembled in
the external level. In sum, the following numbers
of GCCP-pieces are present:

{\bf 1.} 12 GCCP-pieces for each of the types $1a,$ $3b$, $5c$; (subtotal: 36 pieces);

{\bf 2.} \hspace*{1mm} 8 GCCP-pieces for each of the types $2a,4c$; (subtotal: 16 pieces);

{\bf 3.} \hspace*{1mm} 4 GCCP-pieces for each of the types $2b,4b$; (subtotal: 8 pieces).

\noindent Recall that there are 20 cut sections for each of types $wx$, ($x=a,b,c$), and 12 cut sections for each of types $wx$, ($w=[5]$).

Each GCCP-piece is represented in Figure 2 by a corresponding s3p $P$, where our present notation $wa$, resp. $wb$, (resp. $wc$), for $w\in [5]$ is expressed as: either $w_1$ or $w_4$, resp. either $w_2$ or $w_5$, (resp. either $w_3$ or $w_6$), with $w_i$ antipodal in $S$ to $w_{i+3}$, for $i=1,2,3$.
Each such s3p $P$  has endvertices $w_i$ and $w_j$ with a common $w$ but $i\ne j$ in $[6]$. 
Recall that the thick-edge 4-cycles in $\cal G$ containing vertices $1_j$ with $j\neq 2$ are obtained by identification of:

{\bf(i)} the 4 corners bearing symbol $1_5$;

{\bf(ii)} each pair of vertices equally denoted with $z_1^1,y_0^0,y_4^0,x_3^1,x_1^1,y_3^1,y_2^0,x_4^1$;

{\bf(iii)} the s3p pairs with thick edges incident to vertices in item (ii).

\noindent In the interior of each of these 4-cycles there is a sole vertex $w_i$. Let $\cal H$ be the graph whose vertices are all the $w_i$ and such that each edge $w_iw_j$ of $\cal H$ is the s3p $P$ of $\cal G$ with $w_i$ and $w_j$ of degree one in $P$ ($w\in [5]; i,j\in [6]$). This $\cal H$ is representable as the external skeleton of any assemblage of the GCCP and reflects the fact that the 15 GCs split into 5 collections of three pairwise orthogonal GCs each, namely:
$\{\{w_1w_2,w_2w_4,w_4w_5,w_5w_1\}$, $\{w_1w_3,w_3w_4,w_4w_6,w_6w_1\}$, $\{w_2w_3,w_3w_5,
w_5w_6,$ $w_6w_2\}\}$.
Here, each of these three 4-sets of s3ps represents an annulus delimited externally by a GC and formed by 4
corresponding GCCP-pieces, yielding a partial total of 12 GCCP-pieces.
Combining with the notation of Section~\ref{s4}, $G'$ is seen to be induced in $\cal G$ by the thick edges of the
s3ps $P$. Each s3p $P$ of $\cal G$ has an inner vertex of
degree 5 and an inner vertex of degree 3, both in $G'$. Each inner vertex
of an s3p is at distance one in $\cal G$ from just one vertex $w_i$
and at distance 2 in $\cal G$ from just one other vertex $w_j$.
Since $G'$ is a subgraph of $\cal G$, then it is embedded in $S$. The resulting faces of
$G'$ in $S$ are delimited by thick-edge 4-cycles, each containing in its
interior a degree-4 vertex $w_i$ of $\cal G$ that is not in $G'$. Clearly, there is a bijection between those 4-cycles and the vertices $w_i$ contained in their interiors.
\qfd

\end{document}